\documentclass{amsart}
\usepackage{amssymb}

\theoremstyle{definition}

\theoremstyle{remark}

\numberwithin{equation}{section}

\input{tcilatex}

\begin{document}
\title{Linear and nonlinear convolution operator equations on the infinite
strip}
\author{Rishad Shahmurov}
\address{Department of Mathematics, Yeditepe University, Kayishdagi Caddesi,
34755 Kayishdagi, Istanbul, Turkey}
\email{shahmurov@hotmail.com}
\subjclass[2000]{Primary 45N05, 47D06, 35J70 }
\date{}
\keywords{Banach--valued $L_{p}$ spaces, operator-valued multipliers, UMD
spaces, convolution operator equation, differential--operator equation}

\begin{abstract}
In the present paper we characterize the existence and uniqueness of maximal 
$L_{p}$-regular solutions of high order convolution operator equations.
Particularly, we get coercive uniform estimates with respect to spectral
parameter and we show that corresponding realization operator is $R$%
-positive and generates analytic semigroup in $L_{p}$. Then we apply these
results to various problems of nonlinear integro-differential equations.
\end{abstract}

\maketitle



%    General info

\section*{1. Introduction}

The aim of this study is to obtain global existence results for wide class
of convolution operator equations (COE) 
\begin{equation*}
\frac{\partial u}{\partial t}+\sum\limits_{k=0}^{l}\left( a_{k}\ast \frac{%
\partial ^{k}u}{\partial x^{k}}+b_{k}\frac{\partial ^{k}u}{\partial x^{k}}%
\right) +\mu \ast Au+\nu Au=F\left( u,\frac{\partial u}{\partial x},\cdot
\cdot \cdot ,\frac{\partial ^{i}u}{\partial x^{i}}\right) ,\eqno(1.1)
\end{equation*}%
\begin{equation*}
u(0,x)=u_{0}\text{, }t\in \left( 0,T\right) ,\text{ }x\in \left( -\infty
,\infty \right) ,
\end{equation*}%
and%
\begin{equation*}
-\frac{\partial ^{2}u}{\partial t^{2}}+\sum\limits_{k=0}^{l}\left( a_{k}\ast 
\frac{\partial ^{k}u}{\partial x^{k}}+b_{k}\frac{\partial ^{k}u}{\partial
x^{k}}\right) +\mu \ast Au+\nu Au=F\left( u,\frac{\partial u}{\partial t}%
\right) ,\eqno(1.2)
\end{equation*}%
\begin{equation*}
\alpha _{1}u(0,x)+\beta _{1}\frac{\partial u(0,x)}{\partial t}=f_{1}(x)\text{%
, }\alpha _{2}u(T,x)+\beta _{2}\frac{\partial u(T,x)}{\partial t}=f_{2}(x),
\end{equation*}%
where $A$ is a linear operator in a $UMD$ space $E$, $i<l,$ $b_{k},\nu \in 
\mathbf{C}$ and $a_{k}=a_{k}(x),$ $\mu =\mu \left( x\right) \in S^{\prime
}(R,\mathbf{C})$. Equations of above type arise, for instance in the study
of heat flow in materials of fading memory type as well as some equations of
population dynamics. For detailed information on this subject see [1-5].

First we investigate corresponding linear convolution equation with spectral
parameter i.e. 
\begin{equation*}
\left( L+\lambda \right) u=\sum\limits_{k=0}^{l}\left( a_{k}\ast \frac{%
\partial ^{k}u}{\partial x^{k}}+b_{k}\frac{\partial ^{k}u}{\partial x^{k}}%
\right) +\mu \ast Au+\nu Au+\lambda u=f,\eqno(1.3)
\end{equation*}%
in $L_{p}(R;E).$ Particularly, we get coercive uniform $L_{p}$-estimates
with respect to $\lambda ,$ which in its turn implies $R$-positivity of $L.$
The main tool we implement here is an operator-valued Fourier multiplier
theorem (FMT) in $L_{p}(R;E)$. The exposition of FMT, their applications and
some related references can be found in $\left[ 6-9\right] $. For the
references concerning FMT in periodic function spaces, optimal regularity
results for convolution operator equations (COE) and delay differential
operator equations see e.g. [10-14] and reference therein.

In the last section we utilize abstract global and local existence results
in [15-16] to handle nonlinear problems (1.1) and (1.2).

Let $\alpha =(\alpha _{1},\alpha _{2},\cdots ,\alpha _{n}),$ where $\alpha
_{i}$ are integers. An $E$--valued generalized function $D^{\alpha }f$ is
called a generalized derivative in the sense of Schwartz distributions, if
the equality 
\begin{equation*}
<D^{\alpha }f,\varphi >=(-1)^{|\alpha |}<f,D^{\alpha }\varphi >
\end{equation*}%
holds for all $\varphi \in S.$ We indicate mixed derivative as follows 
\begin{equation*}
D^{\alpha }=D_{1}^{\alpha _{1}}D_{2}^{\alpha _{2}}\cdots D_{n}^{\alpha
_{n}},~D_{k}^{i}=\left( \frac{\partial }{\partial x_{k}}\right) ^{i}.
\end{equation*}

The Fourier transform $F:S(X)\rightarrow S(X)$ defined by 
\begin{equation*}
(Ff)(t)\equiv \hat{f}(t)=\int\limits_{R^{N}}\exp (-its)f(s)ds
\end{equation*}%
is an isomorphism whose inverse is given by 
\begin{equation*}
(F^{-1}f)(t)\equiv \check{f}(t)=(2\pi )^{-N}\int\limits_{R^{N}}\exp
(its)f(s)ds,
\end{equation*}%
where $f\in S(X)$ and $t\in R^{N}.$ It is clear that 
\begin{equation*}
F(D_{x}^{\alpha }f)=(i\xi _{1})^{\alpha _{1}}\cdots (i\xi _{n})^{\alpha _{n}}%
\hat{f},~D_{\xi }^{\alpha }(F(f))=F[(-ix_{n})^{\alpha _{1}}\cdots
(-ix_{n})^{\alpha _{n}}f]
\end{equation*}%
for all $f\in S^{\dagger }(R^{n};E).$

Suppose $\Omega \subset R$ and $E_{1}\hookrightarrow E_{0},$ where $%
\hookrightarrow $ denotes continuous and dense injection$.$ $W_{p}^{l}\left(
\Omega ;E_{1},E_{0}\right) $ is a space of all functions $u\in L_{p}\left(
\Omega ;E_{0}\right) $ such that $u^{\left( k\right) }\in L_{p}\left( \Omega
;E\right) ,$ $k=0,\cdot \cdot \cdot ,l$ and 
\begin{equation*}
\left\Vert u\right\Vert _{W_{p}^{l}\left( \Omega ;E_{1},E_{0}\right)
}=\left\Vert u\right\Vert _{L_{p}\left( \Omega ;E_{1}\right)
}+\sum\limits_{k=0}^{l}\left\Vert u^{\left( k\right) }\right\Vert
_{L_{p}\left( \Omega ;E_{0}\right) }<\infty .
\end{equation*}%
For $E_{0}=E_{1}$ the space $W_{p}^{l}\left( \Omega ;E_{1},E_{0}\right) $
will be denoted by $W_{p}^{l}\left( \Omega ;E_{0}\right) .$

Suppose $E_{1}$ and\ $E_{2}$ are two Banach spaces. $B\left(
E_{1},E_{2}\right) $ will denote the space of all bounded linear operators
from $E_{1}$ to $E_{2}.$

A linear operator\ $A$ is said to be $\varphi $-positive in\ a Banach\ space 
$E$,\ with bound $M$ if\ $D\left( A\right) $ is dense in $E$ and 
\begin{equation*}
\left\Vert \left( A+\lambda I\right) ^{-1}\right\Vert _{B\left( E\right)
}\leq M\left( 1+\left\vert \lambda \right\vert \right) ^{-1}
\end{equation*}%
for all $\lambda \in S_{\varphi },$ with $\varphi \in \left[ 0,\left. \pi
\right) \right. ,$where $M$ is a positive constant and $I$ is identity
operator in $E.$

$E\left( A^{\theta }\right) $ denotes the space $D\left( A^{\theta }\right) $
with graphical norm 
\begin{equation*}
\left\Vert u\right\Vert _{E\left( A^{\theta }\right) }=\left( \left\Vert
u\right\Vert ^{p}+\left\Vert A^{\theta }u\right\Vert ^{p}\right) ^{\frac{1}{p%
}},1\leq p<\infty ,-\infty <\theta <\infty .
\end{equation*}%
\vspace{3mm}

\section*{2. Basic notations and parameter depended FMT}

In the next section, we will study COE with spectral parameter via FMT in $%
L_{p}$ spaces. Since we will deal with family of uniformly bounded
functions, we need to formulate slightly extended version of FMT due to
Strkalj and Weis [10]. Let us first introduce some basic definitions and
facts.

\vspace{3mm}

\textbf{Definition 2.0}$.$\textbf{\ }A Banach space $X$ is called $UMD$
space if $X$-valued martingale difference sequences are unconditional in $%
L_{p}\left( R^{d};X\right) $ for $p\in \left( 1,\infty \right) .$ i.e. there
exists a positive constant $C_{p}$ such that for any martingale $\left\{
f_{k},k\in \mathbf{N}_{0}\right\} $, any choice of signs $\left\{
\varepsilon _{k},k\in \mathbf{N}\right\} \in \left\{ -1,1\right\} $ and $%
N\in \mathbf{N}$%
\begin{equation*}
\left\Vert f_{0}+\sum\limits_{k=1}^{N}\varepsilon _{k}\left(
f_{k}-f_{k-1}\right) \right\Vert _{L_{p}\left( \Omega ,\Sigma ,\mu ,X\right)
}\leq C_{p}\left\Vert f_{N}\right\Vert _{L_{p}\left( \Omega ,\Sigma ,\mu
,X\right) }.
\end{equation*}%
It is shown in [17] and $\left[ 18\right] $ that\ a Hilbert operator 
\begin{equation*}
\left( Hf\right) \left( x\right) =\lim\limits_{\varepsilon \rightarrow
0}\int\limits_{\left\vert x-y\right\vert >\varepsilon }\frac{f\left(
y\right) }{x-y}dy
\end{equation*}%
is bounded in the space $L_{p}\left( R,X\right) ,$ $p\in \left( 1,\infty
\right) $ for only those spaces $X,$ which possess the $UMD$ property. $UMD$
spaces include e.g. $L_{p}$, $l_{p}$ spaces and Lorentz spaces $L_{pq},$ $p,$
$q\in \left( 1,\infty \right) $.

\vspace{3mm}

\textbf{Definition 2.1}$.$ \bigskip Let $X$ and $Y$ be Banach spaces. A
family of operators $\tau \subset B\left( X,Y\right) $ is called $R$-bounded
(see e.g. $\left[ 6\right] $ ) if there is a positive constant $C$ and $p\in
\lbrack 1,\infty )$ such that for each $N\in \mathbf{N,}$ $T_{j}\in \tau $, $%
x_{j}\in X$ and for all independent, symmetric, $\left\{ -1,1\right\} -$%
valued random variables $r_{j}$ on a probability space $(\Omega ,\Sigma ,\mu
)$ the inequality%
\begin{equation*}
\left\Vert \sum\limits_{j=1}^{N}r_{j}T_{j}x_{j}\right\Vert _{L_{p}(\Omega
,Y)}\leq C\left\Vert \sum\limits_{j=1}^{N}r_{j}x_{j}\right\Vert
_{L_{p}(\Omega ,X)},
\end{equation*}%
is valid. The smallest such $C$ is called $R$-bound of $\tau ,$ we denote it
by $R_{p}(\tau ).$

The basic properties of $R$-boundedness are collected in the recent
monograph of Denk et al. [6]. For the reader's convenience, we present some
results from [6].

\vspace{3mm}

(a)The definition of $R$-boundedness is independent of $p\in \lbrack
1,\infty ).$

(b) If $\tau \subset B\left( X,Y\right) $ is $R$-bounded then it is
uniformly bounded with%
\begin{equation*}
\sup \left\{ \left\Vert T\right\Vert :\text{ }T\in \tau \right\} \leq
R_{p}(\tau ).
\end{equation*}

(c) If $X$ and $Y$ are Hilbert spaces, $\tau \subset B\left( X,Y\right) $ is 
$R$-bounded $\Longleftrightarrow \tau $ is uniformly bounded.

(d) Let $X,$ $Y$ be Banach spaces and $\tau _{1}$, $\tau _{2}\subset B(X,Y)$
be $R$-bounded. Then%
\begin{equation*}
\tau _{1}+\tau _{2}=\left\{ T+S:T\in \tau _{1},\text{ }S\in \tau _{2}\right\}
\end{equation*}%
is $R$-bounded as well, and $R_{p}(\tau _{1}+\tau _{2})\leq R_{p}(\tau
_{1})+R_{p}(\tau _{2})$ .

(e) Let $X,$ $Y,$ $Z$ be Banach spaces and $\tau _{1}\subset B(X,Y)$ and $%
\tau _{2}\subset B(Y,Z)$ be $R$-bounded. Then%
\begin{equation*}
\tau _{1}\tau _{2}=\left\{ ST:T\in \tau _{1},\text{ }S\in \tau _{2}\right\}
\end{equation*}%
is $R$-bounded as well, and $R_{p}(\tau _{1}\tau _{2})\leq R_{p}(\tau
_{1})R_{p}(\tau _{2})$.

One of the most important tools in $R$-boundedness is the contraction
principle of Kahane. We shall frequently apply it in the next sections.

\vspace{3mm}

\textbf{[6}, \textbf{Lemma 3.5.] }Let $X$ be a Banach spaces, $n\in N,$ $%
x_{j}\in X,$ $r_{j}$ independent, symmetric, $\left\{ -1,1\right\} $-valued
random variables on a probability space $(\Omega ,\Sigma ,\mu )$ and $\alpha
_{j,}\beta _{j}\in \mathbf{C}$ such that $\left\vert \alpha _{j}\right\vert
\leq \left\vert \beta _{j}\right\vert ,$ for each $j=1,\cdot \cdot \cdot ,N.$
Then%
\begin{equation*}
\left\Vert \sum\limits_{j=1}^{N}\alpha _{j}r_{j}x_{j}\right\Vert
_{L_{p}(\Omega ,X)}\leq 2\left\Vert \sum\limits_{j=1}^{N}\beta
_{j}r_{j}x_{j}\right\Vert _{L_{p}(\Omega ,X)}.
\end{equation*}%
The constant 2 can be omitted in case where $\alpha _{j}$ and $\beta _{j}$
are real.

\vspace{3mm}

\textbf{Definition 2.2. }A family of uniformly bounded functions $%
m_{h}:R^{d}\rightarrow \mathbf{C}$ is called a uniform collection of Fourier
multipliers (UFM) if there exists a positive constant $C>0,$ independent of
parameter $h\in Q$, such that

\begin{equation*}
\left\Vert F^{-1}\left[ m_{h}\hat{f}\right] \right\Vert _{L_{p}\left(
R^{d};X\right) }\leq C\left\Vert f\right\Vert _{L_{p}\left( R^{d};X\right) }%
\eqno(2.1)\ 
\end{equation*}%
for all $f\in S\left( R^{d};X\right) .$

The set of all $L_{p}(X)$-UFM will be denoted by $M_{p}(X)$ and the smallest
constant $C$ satisfying (2.1) by $\left\Vert m_{h}\right\Vert _{M_{p}(X)}.$

\vspace{3mm}

\textbf{Theorem 2.3. }Let $X$ be a $UMD$ space$.$ Then for any $p\in \left(
1,\infty \right) $ there is a constant $C<\infty $ such that for all family
of uniformly bounded functions $m_{h}:R^{d}\rightarrow \mathbf{C}$ whose
distributional derivatives $D^{\alpha }m_{h}$ of order $\alpha \leq
(1,...,1) $ are represented by functions, we have%
\begin{equation*}
\left\Vert m_{h}\right\Vert _{M_{p}(X)}\leq C\sup_{h}\sup \left\{ \left\vert
\xi \right\vert ^{\left\vert \alpha \right\vert }\left\vert D^{\alpha
}m_{h}(\xi )\right\vert ,\text{ }\xi \in R^{d}\backslash \left\{ 0\right\} ,%
\text{ }\alpha \leq (1,...,1)\right\} .
\end{equation*}

It is easy to see that Theorem 2.3 is a trivial consequence of [9,
Proposition 2]. Since in [9, Proposition 2] author finds a constant $C$
independent of the multiplier functions, we can take supremum over
parameters to obtain above theorem.

The following result is a parameter depended version of operator valued
Miklin's theorem [10, Theorem 4.4].

\vspace{3mm}

\textbf{Theorem 2.4. }Let $X$ and $Y$ be $UMD$ spaces and $1<p<\infty .$ If
the family of operator-valued function $M_{h}:R^{d}\backslash \left\{
0\right\} \rightarrow B(X,Y)$ has the property that their distributional
derivatives $D^{\alpha }M_{h}$ of order $\alpha \leq (1,...,1)$ are
represented by functions and 
\begin{equation*}
\sup_{h}R\left\{ \left\vert x\right\vert ^{\left\vert \alpha \right\vert
}D^{\alpha }M_{h}(x);\text{ }x\in R^{d}\backslash \left\{ 0\right\} ,\text{ }%
(\alpha \leq 1,...,1)\text{ }\right\} <\infty \text{ }
\end{equation*}%
holds, then $M_{h}$ $\ $is UFM.

Taking into consideration the estimate (4.4) of [10, Theorem 4.4] and by
using the similar reasoning as in Theorem 2.3 we get the assertion of
Theorem 2.4.

\vspace{3mm}

\section*{3. Maximal regularity for (1.3)}

Now let us consider the high order integro-differential equation (1.3) in $%
L_{p}(R;E).$ Here we derive a sufficient condition which guarantee the
maximal regularity of (1.3).

\vspace{3mm}

\textbf{Definition 3.0. }Let $E$ be a Banach space and $D\left( A\right) $
dense in $E.$ A $\varphi $-positive operator $A$ is said to be $R$-positive
if the following set%
\begin{equation*}
\left\{ (1+\xi )(A+\xi )^{-1}:\text{ }\xi \in S_{\varphi }\right\}
\end{equation*}%
is $R$-bounded.

In what follows $\dot{R}$ will denote the set of real numbers excluding zero
i.e. $\dot{R}=R/\left\{ 0\right\} .$

\bigskip

\textbf{Condition 3.1. }Let $E$ be a Banach space and $A$ be a $R$-positive
operator in $E.$ Suppose the following are satisfied:

(\textbf{1}) $b_{k},\nu \in \mathbf{C,}$ $a_{k},\mu \in S^{\prime }(R,%
\mathbf{C}),$ $\hat{a}_{k},\hat{\mu}\in C^{1}(\dot{R},\mathbf{C})$ and 
\begin{equation*}
C_{\mu }=\inf_{\xi \in \dot{R}}\left\vert \hat{\mu}(\xi )+\nu \right\vert >0%
\text{ };
\end{equation*}

(\textbf{2}) There exist a constant $C_{N}$ such that%
\begin{equation*}
\left\vert N\left( \xi \right) \right\vert =\left\vert
\sum\limits_{k=0}^{l}\left( b_{k}+\hat{a}_{k}(\xi )\right) (i\xi
)^{k}\right\vert \geq C_{N}\left\vert \xi \right\vert ^{l};
\end{equation*}

(\textbf{3}) 
\begin{equation*}
\eta \left( \xi \right) =\frac{N\left( \xi \right) }{\hat{\mu}(\xi )+\nu }%
\in S_{\varphi _{1}},\text{ }\varphi _{1}\in \left[ 0,\right. \left. \pi
\right) \text{ and }\lambda \in S_{\varphi _{2}},\varphi _{2}\in \lbrack
0,\left. \pi \right)
\end{equation*}%
are so that $\varphi _{1}+\varphi _{2}<\pi $ and 
\begin{equation*}
\lambda +\eta \left( \xi \right) \in S_{\varphi },
\end{equation*}

(\textbf{4}) There are constants $C_{1}$ and $C_{2}$ so that 
\begin{equation*}
\left\vert \xi ^{m}\frac{d^{m}}{d\xi ^{m}}\hat{a}_{k}(\xi )\right\vert \leq
C_{1},\text{ for all }k=0,1,\cdot \cdot \cdot ,l,
\end{equation*}%
\begin{equation*}
\left\vert \text{ }\xi ^{m}\frac{d^{m}}{d\xi ^{m}}\hat{\mu}(\xi )\right\vert
\leq C_{2},\text{ }
\end{equation*}%
for $m=0,1$ and $\xi \in \dot{R}.$

Taking into account Condition 3.1 and by using the Kahane's contraction
principle we shall estimate $R$-bounds of the following sets:%
\begin{equation*}
\left\{ m_{i}\left( \xi ,\lambda \right) :\xi \in \dot{R}\right\}
\end{equation*}%
and 
\begin{equation*}
\left\{ \xi \frac{d}{d\xi }m_{i}\left( \xi ,\lambda \right) :\xi \in \dot{R}%
\right\} ,\text{ }i=0,\cdot \cdot \cdot ,4,
\end{equation*}%
where%
\begin{equation*}
m_{0}(\xi ,\lambda )=(\hat{\mu}(\xi )+\nu )^{-1}\left[ A+\eta \left( \xi
\right) +\lambda )\right] ^{-1},
\end{equation*}%
\begin{eqnarray*}
m_{1}\left( \xi ,\lambda \right) &=&\sum\limits_{k=0}^{l}\left\vert \lambda
\right\vert ^{1-\frac{k}{l}}\left( i\xi \right) ^{k}m_{0}(\xi ,\lambda ),%
\text{ }m_{2}\left( \xi ,\lambda \right) =Am_{0}(\xi ,\lambda ),\text{ } \\
m_{3}\left( \xi ,\lambda \right) &=&\sum\limits_{k=0}^{l}\left\vert \lambda
\right\vert ^{1-\frac{k}{l}}\hat{a}_{k}(\xi )\left( i\xi \right)
^{k}m_{0}(\xi ,\lambda )\text{ and }m_{4}\left( \xi ,\lambda \right) =\hat{%
\mu}(\xi )Am_{0}(\xi ,\lambda )
\end{eqnarray*}%
The main theorem of this section is based on the following lemmas.

\vspace{3mm}

\textbf{Lemma 3.2. }Assume there are some constants $C_{1}$ and $C_{2}$ so
that 
\begin{equation*}
\left\vert \hat{a}_{k}(\xi )\right\vert \leq C_{1}\text{ for all }k=0,\cdot
\cdot \cdot ,l
\end{equation*}%
and%
\begin{equation*}
\text{ }\left\vert \hat{\mu}(\xi )\right\vert \leq C_{2}
\end{equation*}%
for all $\xi \in \dot{R}.$ If (1)-(3) of Condition 3.1 hold then $\left\{
m_{i}\left( \xi ,\lambda \right) :\xi \in \dot{R}\right\} ,$ $i=0,\cdot
\cdot \cdot ,4$ are $R$-bounded.

\textbf{Proof. }Since $\eta \left( \xi \right) \in S_{\varphi _{1}}$, $%
\varphi _{1}\in \left[ 0,\right. \left. \pi \right) $ and $\lambda \in
S_{\varphi _{2}}$ for $\varphi _{2}\in \left[ 0,\right. \left. \pi \right) ,$%
\ by [19, Lemma 2.3] there exist $K>0$ independent of $\xi $ so that 
\begin{equation*}
\left\vert 1+\eta \left( \xi \right) +\lambda \right\vert ^{-1}\leq
K(1+\left\vert \eta \left( \xi \right) \right\vert +\left\vert \lambda
\right\vert )^{-1}.
\end{equation*}%
Therefore, for all $\xi \in R$ we have uniform estimate%
\begin{equation*}
\begin{array}{lll}
\frac{1}{\left\vert \hat{\mu}(\xi )+\nu \right\vert }\times \frac{1}{%
\left\vert 1+\eta \left( \xi \right) +\lambda \right\vert } & \leq & %
\displaystyle K\frac{1}{\left\vert \hat{\mu}(\xi )+\nu \right\vert
+\left\vert \lambda \right\vert \left\vert \hat{\mu}(\xi )+\nu \right\vert
+\left\vert N\left( \xi \right) \right\vert } \\ 
&  &  \\ 
& \leq & \displaystyle K\frac{1}{C_{\mu }+\left\vert \lambda \right\vert
C_{\mu }+C_{N}\left\vert \xi \right\vert ^{l}}\leq \frac{K}{C_{\mu }}.%
\end{array}%
\end{equation*}%
Now let us define families of operators 
\begin{equation*}
\tau =\left\{ T_{j}=(1+\eta \left( \xi _{j}\right) +\lambda )\left( A+\eta
\left( \xi _{j}\right) +\lambda \right) ^{-1};\text{ }\xi _{j}\in \dot{R}%
\right\} ,
\end{equation*}%
and%
\begin{equation*}
\tau _{i}=\left\{ T_{j}^{i}=m_{i}(\xi _{j},\lambda );\text{ }\xi _{j}\in 
\dot{R}\right\} ,\text{ }i=0,\cdot \cdot \cdot ,4.
\end{equation*}%
Taking into consideration $R$-positivity of $A$, applying the assumptions of
Condition 3.1 and the Kahane's contraction principle we get desired result:%
\begin{equation*}
\begin{array}{lll}
\left\Vert \sum\limits_{j=1}^{N}r_{j}T_{j}^{0}x_{j}\right\Vert _{X} & = & %
\displaystyle\left\Vert \sum\limits_{j=1}^{N}r_{j}\frac{1}{\left( \hat{\mu}%
(\xi _{j})+\nu \right) \left( 1+\eta \left( \xi _{j}\right) +\lambda \right) 
}T_{j}x_{j}\right\Vert _{X} \\ 
&  &  \\ 
& \leq & \displaystyle2\frac{K}{C_{\mu }}\left\Vert
\sum\limits_{j=1}^{N}r_{j}T_{j}x_{j}\right\Vert _{X}\leq 2\frac{K}{C_{\mu }}%
R_{p}(\tau )\left\Vert \sum\limits_{j=1}^{N}r_{j}x_{j}\right\Vert _{X},%
\end{array}%
\end{equation*}%
where $X=L_{p}((0,1),E).$ Hence 
\begin{equation*}
R_{p}(\tau _{0})\leq 2\frac{K}{C_{\mu }}R_{p}(\tau ).
\end{equation*}%
It is clear that 
\begin{equation*}
\sum\limits_{k=0}^{l}\left\vert y\right\vert ^{k}\leq l(1+\left\vert
y\right\vert ^{l})
\end{equation*}%
for any $y\in \mathbf{C}.$ Thus making use of the above inequality we get%
\begin{equation*}
\left\vert \sum\limits_{k=0}^{l}\left\vert \lambda \right\vert ^{1-\frac{k}{l%
}}\left( i\xi \right) ^{k}\right\vert \leq \left\vert \lambda \right\vert
\sum\limits_{k=0}^{l}\left( \left\vert \lambda \right\vert ^{-\frac{1}{l}%
}\left\vert \xi \right\vert \right) ^{k}\leq l\left( \left\vert \lambda
\right\vert +\left\vert \xi \right\vert ^{l}\right) ,
\end{equation*}%
\begin{equation*}
\begin{array}{lll}
\frac{1}{\left\vert \hat{\mu}(\xi )+\nu \right\vert }\frac{\left\vert
\sum\limits_{k=0}^{l}\left\vert \lambda \right\vert ^{1-\frac{k}{l}}\left(
i\xi \right) ^{k}\right\vert }{\left\vert 1+\eta \left( \xi \right) +\lambda
\right\vert } & \leq & \displaystyle Kl\frac{\left\vert \lambda \right\vert
+\left\vert \xi \right\vert ^{l}}{\left\vert \hat{\mu}(\xi )+\nu \right\vert
+\left\vert \lambda \right\vert \left\vert \hat{\mu}(\xi )+\nu \right\vert
+\left\vert N\left( \xi \right) \right\vert } \\ 
&  &  \\ 
& \leq & \displaystyle Kl\frac{\left\vert \lambda \right\vert +\left\vert
\xi \right\vert ^{l}}{C_{\mu }+\left\vert \lambda \right\vert C_{\mu
}+C_{N}\left\vert \xi \right\vert ^{l}}\leq KlM%
\end{array}%
\end{equation*}%
and%
\begin{equation*}
\frac{1}{\left\vert \hat{\mu}(\xi )+\nu \right\vert }\frac{\left\vert
\sum\limits_{k=0}^{l}\left\vert \lambda \right\vert ^{1-\frac{k}{l}}\hat{a}%
_{k}(\xi )\left( i\xi \right) ^{k}\right\vert }{\left\vert 1+\eta \left( \xi
\right) +\lambda \right\vert }\leq KlMC_{1},
\end{equation*}%
where $M^{-1}=\min \left\{ C_{\mu },C_{N}\right\} .$ Again by the Kahane's
contraction principle we have%
\begin{equation*}
\begin{array}{lll}
\left\Vert \sum\limits_{j=1}^{N}r_{j}T_{j}^{1}x_{j}\right\Vert _{X} & = & %
\displaystyle\left\Vert \sum\limits_{j=1}^{N}r_{j}\frac{\sum%
\limits_{k=0}^{l}\left\vert \lambda \right\vert ^{1-\frac{k}{l}}\left( i\xi
_{j}\right) ^{k}}{\left( \hat{\mu}(\xi _{j})+\nu \right) \left( 1+\eta
\left( \xi _{j}\right) +\lambda \right) }T_{j}x_{j}\right\Vert _{X} \\ 
&  &  \\ 
& \leq & \displaystyle2KlM\left\Vert
\sum\limits_{j=1}^{N}r_{j}T_{j}x_{j}\right\Vert _{X}\leq 2KlMR_{p}(\tau
)\left\Vert \sum\limits_{j=1}^{N}r_{j}x_{j}\right\Vert _{X},%
\end{array}%
\end{equation*}%
and%
\begin{equation*}
\begin{array}{lll}
\left\Vert \sum\limits_{j=1}^{N}r_{j}T_{j}^{3}x_{j}\right\Vert _{X} & = & %
\displaystyle\left\Vert \sum\limits_{j=1}^{N}r_{j}\frac{\sum%
\limits_{k=0}^{l}\left\vert \lambda \right\vert ^{1-\frac{k}{l}}\hat{a}%
_{k}\left( i\xi _{j}\right) ^{k}}{\left( \hat{\mu}(\xi _{j})+\nu \right)
\left( 1+\eta \left( \xi _{j}\right) +\lambda \right) }T_{j}^{0}x_{j}\right%
\Vert _{X} \\ 
&  &  \\ 
& \leq & \displaystyle2KlMC_{1}R_{p}(\tau )\left\Vert
\sum\limits_{j=1}^{N}r_{j}x_{j}\right\Vert _{X},%
\end{array}%
\end{equation*}%
which implies%
\begin{equation*}
R_{p}(\tau _{1})\leq 2KlMR_{p}(\tau )
\end{equation*}%
and%
\begin{equation*}
R_{p}(\tau _{3})\leq 2KlMC_{1}R_{p}(\tau ).
\end{equation*}%
Finally, by virtue of the resolvent property we get 
\begin{equation*}
\begin{array}{lll}
R_{p}(\tau _{2}) & \leq & \displaystyle R_{p}\left\{ I\frac{1}{\hat{\mu}(\xi
_{j})+\nu };\text{ }\xi _{j}\in \dot{R}\right\} \\ 
&  &  \\ 
& + & \displaystyle R_{p}\left\{ \left( \eta \left( \xi _{j}\right) +\lambda
\right) \left( A+\eta \left( \xi _{j}\right) +\lambda \right) ^{-1};\text{ }%
\xi _{j}\in \dot{R}\right\} \\ 
&  &  \\ 
& \leq & \displaystyle2\left( \frac{1}{C_{\mu }}+R_{p}(\tau )\right)%
\end{array}%
\end{equation*}%
and 
\begin{equation*}
\begin{array}{lll}
R_{p}(\tau _{4}) & \leq & \displaystyle2C_{2}\left( \frac{1}{C_{\mu }}%
+R_{p}(\tau )\right) .%
\end{array}%
\end{equation*}%
Next we estimate the derivatives of operator valued functions $m_{i}(\xi ).$

\vspace{3mm}

\textbf{Lemma 3.3. }If the Condition 3.1 holds then $\left\{ \xi \frac{d}{%
d\xi }m_{i}\left( \xi ,\lambda \right) :\xi \in \dot{R}\right\} ,$ $%
i=0,\cdot \cdot \cdot ,4$ are $R$-bounded\textbf{.}%
\begin{equation*}
\end{equation*}%
\textbf{Proof. }For the sake of simplicity we only estimate $\left\{ \xi 
\frac{d}{d\xi }m_{1}\left( \xi ,\lambda \right) :\xi \in \dot{R}\right\} $%
and $\left\{ \xi \frac{d}{d\xi }m_{2}\left( \xi ,\lambda \right) :\xi \in 
\dot{R}\right\} $. Computing the first derivative of $m_{1}$ we obtain 
\begin{equation*}
\begin{array}{lll}
R_{p}\left( \left\{ \xi \frac{d}{d\xi }m_{1}\left( \xi ,\lambda \right) :\xi
\in \dot{R}\right\} \right) & \leq & \displaystyle\dsum%
\limits_{i=1}^{5}R_{p}\left( \tau _{i}\right) ,%
\end{array}%
\end{equation*}%
where 
\begin{equation*}
\begin{array}{lll}
\tau _{1} & = & \displaystyle\left\{ \sum\limits_{k=1}^{l}\left\vert \lambda
\right\vert ^{1-\frac{k}{l}}k\left( i\xi \right) ^{k-1}m_{0}(\xi ,\lambda
):\xi \in \dot{R}\right\} , \\ 
&  &  \\ 
\tau _{2} & = & \displaystyle\left\{ \frac{\frac{d\hat{\mu}}{d\xi }}{\hat{\mu%
}(\xi )+\nu }m_{1}\left( \xi ,\lambda \right) :\xi \in \dot{R}\right\} , \\ 
&  &  \\ 
\tau _{3} & = & \displaystyle\left\{ m_{1}\left( \xi ,\lambda \right)
\sum\limits_{k=0}^{l}\frac{d\hat{a}_{k}}{d\xi }\left( i\xi \right)
^{k}m_{0}(\xi ,\lambda ):\xi \in \dot{R}\right\} , \\ 
&  &  \\ 
\tau _{4} & = & \displaystyle\left\{ m_{1}\left( \xi ,\lambda \right)
\sum\limits_{k=1}^{l}\left( b_{k}+\hat{a}_{k}(\xi )\right) k\left( i\xi
\right) ^{k-1}m_{0}(\xi ,\lambda ):\xi \in \dot{R}\right\}%
\end{array}%
\end{equation*}%
and%
\begin{equation*}
\begin{array}{lll}
\tau _{5} & = & \displaystyle\left\{ m_{1}\left( \xi ,\lambda \right) \frac{%
N(\xi )\frac{d\hat{\mu}}{d\xi }}{\hat{\mu}(\xi )+\nu }m_{0}(\xi ,\lambda
):\xi \in \dot{R}\right\} .%
\end{array}%
\end{equation*}%
Making use of similar arguments as in Lemma 3.2 we get%
\begin{equation*}
\begin{array}{lll}
\frac{1}{\left\vert \hat{\mu}(\xi )+\nu \right\vert }\frac{\left\vert
\sum\limits_{k=1}^{l}\left\vert \lambda \right\vert ^{1-\frac{k}{l}}k\left(
i\xi \right) ^{k-1}\right\vert }{\left\vert 1+\eta \left( \xi \right)
+\lambda \right\vert } & \leq & \displaystyle Kl\frac{\left\vert \lambda
\right\vert ^{1-\frac{1}{l}}+\left\vert \lambda \right\vert ^{1-\frac{2}{l}%
}\left\vert \xi \right\vert \cdot \cdot \cdot +\left\vert \xi \right\vert
^{l-1}}{C_{\mu }+C_{\mu }\left\vert \lambda \right\vert +C_{N}\left\vert \xi
\right\vert ^{l}}\leq KlM_{1}, \\ 
&  &  \\ 
\frac{\left\vert \frac{d\hat{\mu}}{d\xi }\right\vert }{\left\vert \hat{\mu}%
(\xi )+\nu \right\vert } & \leq & \displaystyle\frac{C_{2}}{C_{\mu }}, \\ 
&  &  \\ 
\frac{1}{\left\vert \hat{\mu}(\xi )+\nu \right\vert }\frac{\left\vert
\sum\limits_{k=0}^{l}\frac{d\hat{a}_{k}}{d\xi }\left( i\xi \right)
^{k}\right\vert }{\left\vert 1+\eta \left( \xi \right) +\lambda \right\vert }
& \leq & \displaystyle Kl\frac{C_{1}}{C_{\mu }}, \\ 
&  &  \\ 
\frac{1}{\left\vert \hat{\mu}(\xi )+\nu \right\vert }\frac{\left\vert
\sum\limits_{k=1}^{l}\left( b_{k}+\hat{a}_{k}(\xi )\right) k\left( i\xi
\right) ^{k-1}\right\vert }{\left\vert 1+\eta \left( \xi \right) +\lambda
\right\vert } & \leq & \displaystyle KlM_{1}(C_{b}+C_{1}),%
\end{array}%
\end{equation*}%
and%
\begin{equation*}
\begin{array}{lll}
\left\vert \frac{N(\xi )\frac{d\hat{\mu}}{d\xi }}{\hat{\mu}(\xi )+\nu }%
\right\vert \frac{1}{\left\vert \hat{\mu}(\xi )+\nu \right\vert }\frac{1}{%
\left\vert 1+\eta \left( \xi \right) +\lambda \right\vert } & \leq & %
\displaystyle K\frac{C_{2}}{C_{\mu }^{2}}(C_{b}+C_{1}).%
\end{array}%
\end{equation*}%
Thus by the Kahane's contraction principle we have%
\begin{equation*}
\begin{array}{lll}
R_{p}\left( \tau _{1}\right) & \leq & \displaystyle2KlM_{1}R_{p}(\tau ), \\ 
&  &  \\ 
R_{p}\left( \tau _{2}\right) & \leq & \displaystyle2\frac{C_{2}}{C_{\mu }}%
KlMR_{p}(\tau ), \\ 
&  &  \\ 
R_{p}\left( \tau _{3}\right) & \leq & \displaystyle4\frac{C_{1}}{C_{\mu }}%
M\left( KlR_{p}(\tau )\right) ^{2}, \\ 
&  &  \\ 
R_{p}\left( \tau _{4}\right) & \leq & \displaystyle4MM_{1}(C_{b}+C_{1})%
\left( KlR_{p}(\tau )\right) ^{2},%
\end{array}%
\end{equation*}%
and%
\begin{equation*}
\begin{array}{lll}
R_{p}\left( \tau _{5}\right) & \leq & \displaystyle4\frac{C_{2}}{C_{\mu }^{2}%
}(C_{b}+C_{1})lM\left( KR_{p}(\tau )\right) ^{2}.%
\end{array}%
\end{equation*}%
Similarly, to show $\left\{ \xi \frac{d}{d\xi }m_{2}\left( \xi ,\lambda
\right) :\xi \in R\backslash \left\{ 0\right\} \right\} $ is $R$-bounded, it
suffices to estimate the following terms 
\begin{eqnarray*}
&&\left\{ \frac{\sum\limits_{k=0}^{l}\frac{d\hat{a}_{k}}{d\xi }\left( i\xi
\right) ^{k}}{\hat{\mu}(\xi )+\nu }m_{0}(\xi ,\lambda ):\xi \in \dot{R}%
\right\} , \\
&&\left\{ \frac{\sum\limits_{k=1}^{l}\left( b_{k}+\hat{a}_{k}(\xi )\right)
k\left( i\xi \right) ^{k-1}}{\hat{\mu}(\xi )+\nu }m_{0}(\xi ,\lambda ):\xi
\in \dot{R}\right\} ,
\end{eqnarray*}%
and%
\begin{equation*}
\left\{ \frac{N(\xi )\frac{d\hat{\mu}}{d\xi }}{\left( \hat{\mu}(\xi )+\nu
\right) ^{2}}m_{0}(\xi ,\lambda ):\xi \in \dot{R}\right\} .
\end{equation*}%
By virtue of above techniques and Lemma 3.2 one can easily find $R$-bounds
for these sets. The other cases can be proven analogously.\ Hence proof is
completed.

From the Lemma 3.2, Lemma 3.3 and Theorem 2.4 we get the following result.

\vspace{3mm}

\textbf{Corollary 3.4. }Let $E$ be a $UMD$ space and $A$ be a $R$-positive
operator in $E.$ If the Condition 3.1 holds and $1<p<\infty $ then
operator-valued functions $m_{i}\left( \xi ,\lambda \right) $ are UFM in $%
L_{p}\left( R;E\right) .$

Now, making use of Corollary 3.4 we obtain our main result for this section.

\vspace{3mm}

\textbf{Theorem 3.5.}\ Let $E$ be a $UMD$ space and $A$ be a $R$-positive
operator in $E.$ If$\ $the Condition 3.1 is satisfied then for each $f\in
Y=L_{p}(R;E),$ (1.3) has a unique solution $u\in W_{p}^{l}(R;E(A);E)$ and
the coercive uniform estimate 
\begin{equation*}
\sum\limits_{k=0}^{l}\left\vert \lambda \right\vert ^{1-\frac{k}{l}}\left(
\left\Vert a_{k}\ast \frac{d^{k}u}{dt^{k}}\right\Vert _{Y}+\left\Vert \frac{%
d^{k}u}{dt^{k}}\right\Vert _{Y}\right) +\left\Vert \mu \ast Au\right\Vert
_{Y}+\left\Vert Au\right\Vert _{Y}\leq C\left\Vert f\right\Vert _{Y}\eqno%
(3.1)
\end{equation*}%
holds.

\textbf{Proof}. Applying the Fourier transform to equation $\left(
1.3\right) $\ we obtain 
\begin{equation*}
(\hat{\mu}(\xi )+\nu )\left[ A+\eta \left( \xi \right) +\lambda \right] \hat{%
u}\left( \xi \right) =\hat{f}\left( \xi \right) .
\end{equation*}%
Due to positivity of $A$ and Condition 3.1, $A+\eta \left( \xi \right)
+\lambda $ is invertible in $E$. Thus, solutions of the equation $\left(
1.3\right) $ can be represented in the form 
\begin{equation*}
u\left( x\right) =F^{-1}\left( m_{0}\hat{f}\right) \eqno(3.2)
\end{equation*}%
where 
\begin{equation*}
m_{0}(\xi ,\lambda )=(\hat{\mu}(\xi )+\nu )^{-1}R\left( \xi ,\lambda \right) 
\text{ and }R\left( \xi ,\lambda \right) =\left[ A+\eta \left( \xi \right)
+\lambda )\right] ^{-1}.
\end{equation*}%
Thus, by using $\left( 3.2\right) $ we obtain

\begin{equation*}
\sum\limits_{k=0}^{l}\left\vert \lambda \right\vert ^{1-\frac{k}{l}%
}\left\Vert \frac{d^{k}u}{dt^{k}}\right\Vert _{Y}\equiv \left\Vert F^{-1}%
\left[ m_{1}\left( \xi ,\lambda \right) \hat{f}\right] \right\Vert _{Y},
\end{equation*}

\begin{equation*}
\left\Vert Au\right\Vert _{Y}=\left\Vert F^{-1}\left[ m_{2}\left( \xi
,\lambda \right) \hat{f}\right] \right\Vert _{Y},
\end{equation*}%
\begin{equation*}
\sum\limits_{k=0}^{l}\left\vert \lambda \right\vert ^{1-\frac{k}{l}%
}\left\Vert a_{k}\ast \frac{d^{k}u}{dt^{k}}\right\Vert _{Y}\equiv \left\Vert
F^{-1}\left[ m_{3}\left( \xi ,\lambda \right) \hat{f}\right] \right\Vert
_{Y},
\end{equation*}%
and%
\begin{equation*}
\left\Vert \mu \ast Au\right\Vert _{Y}=\left\Vert F^{-1}\left[ m_{4}\left(
\xi ,\lambda \right) \hat{f}\right] \right\Vert _{Y},
\end{equation*}%
where 
\begin{eqnarray*}
m_{1}\left( \xi ,\lambda \right) &=&\sum\limits_{k=0}^{l}\left\vert \lambda
\right\vert ^{1-\frac{k}{l}}\left( i\xi \right) ^{k}m_{0}(\xi ,\lambda ),%
\text{ }m_{2}\left( \xi ,\lambda \right) =Am_{0}(\xi ,\lambda ),\text{ } \\
m_{3}\left( \xi ,\lambda \right) &=&\sum\limits_{k=0}^{l}\left\vert \lambda
\right\vert ^{1-\frac{k}{l}}\hat{a}_{k}(\xi )\left( i\xi \right)
^{k}m_{0}(\xi ,\lambda )\text{ and }m_{4}\left( \xi ,\lambda \right) =\hat{%
\mu}(\xi )Am_{0}(\xi ,\lambda ).
\end{eqnarray*}%
By the Corollary 3.4 $m_{i}$, $i=1,\cdot \cdot \cdot ,4$ are UFM in $%
L_{p}\left( R;E\right) .$ Hence proof is completed.

Let $L$ be an operator generated by the problem $\left( 1.3\right) $ i. e.\ 
\begin{equation*}
Lu=\sum\limits_{k=0}^{l}a_{k}\ast \frac{d^{k}u}{dt^{k}}+\sum%
\limits_{k=0}^{l}b_{k}\frac{d^{k}u}{dt^{k}}+\mu \ast Au+\nu Au,
\end{equation*}%
where%
\begin{equation*}
D\left( L\right) =W_{p}^{l}\left( R;E\left( A\right) ,E\right)
=W_{p}^{l}\left( R;E\right) \cap L_{p}\left( R;E\left( A\right) \right) .
\end{equation*}%
The estimate (3.1) implies that $L$ is bijective and%
\begin{equation*}
\left\Vert u\right\Vert _{W_{p}^{l}\left( R;E\left( A\right) ,E\right) }\leq
C\left\Vert f\right\Vert _{Y}=C\left\Vert Lu\right\Vert _{Y}.
\end{equation*}%
Therefore $L$ has a continuous inverse 
\begin{equation*}
L^{-1}:L_{p}\left( R,E\right) \rightarrow W_{p}^{l}\left( R;E\left( A\right)
,E\right) .\eqno(3.3)
\end{equation*}%
Under some restrictions on coefficients $a_{k}$ and $\mu ,$ it is also
possible to show $L$ is continuous from $W_{p}^{l}\left( R;E\left( A\right)
,E\right) $ to $L_{p}\left( R,E\right) $. Really,%
\begin{equation*}
\left\Vert Lu\right\Vert _{Y}\leq \sum\limits_{k=0}^{l}\left\Vert a_{k}\ast 
\frac{d^{k}u}{dt^{k}}\right\Vert _{Y}+\sum\limits_{k=0}^{l}b_{k}\left\Vert 
\frac{d^{k}u}{dt^{k}}\right\Vert _{Y}+\left\Vert \mu \ast Au\right\Vert
_{Y}+\left\Vert \nu Au\right\Vert _{Y}.
\end{equation*}%
\ Since $\hat{a}_{k}$ and $\hat{\mu}$ are Fourier multipliers, we have%
\begin{equation*}
\sum\limits_{k=0}^{l}\left\Vert a_{k}\ast \frac{d^{k}u}{dt^{k}}\right\Vert
_{Y}=\sum\limits_{k=0}^{l}\left\Vert F^{-1}\left[ \hat{a}_{k}F\left( \frac{%
d^{k}u}{dt^{k}}\right) \right] \right\Vert _{Y}\leq
C\sum\limits_{k=0}^{l}\left\Vert \frac{d^{k}u}{dt^{k}}\right\Vert _{Y}
\end{equation*}%
and%
\begin{equation*}
\left\Vert \mu \ast Au\right\Vert _{Y}\leq C\left\Vert Au\right\Vert _{Y}.
\end{equation*}%
Therefore 
\begin{equation*}
\left\Vert Lu\right\Vert _{Y}\leq C\left\Vert u\right\Vert _{W_{p}^{l}\left(
R;E\left( A\right) ,E\right) }.\eqno(3.4)
\end{equation*}%
Hence by (3.3) and (3.4) the operator $L$ is an isomorphism from $%
W_{p}^{l}\left( R;E\left( A\right) ,E\right) $ to $L_{p}\left( R,E\right) $
and 
\begin{equation*}
\left\Vert Lu\right\Vert _{L_{p}\left( R,E\right) }\approx \left\Vert
u\right\Vert _{W_{p}^{l}\left( R;E\left( A\right) ,E\right) }
\end{equation*}%
for all $u\in W_{p}^{l}\left( R;E\left( A\right) ,E\right) .$

Now let us write the above observation as a separate corollary.

\vspace{3mm}

\textbf{Result 3.6. }Let $E$ be a $UMD$ space and $A$ be a $R$-positive
operator in $E.$ If Condition 3.1 holds then $L$ is an isomorphism from $%
W_{p}^{l}\left( R;E\left( A\right) ,E\right) $ to $L_{p}\left( R,E\right) $
i.e. 
\begin{equation*}
\left\Vert Lu\right\Vert _{L_{p}\left( R,E\right) }\approx \left\Vert
u\right\Vert _{W_{p}^{l}\left( R;E\left( A\right) ,E\right) }
\end{equation*}%
for all $u\in W_{p}^{l}\left( R;E\left( A\right) ,E\right) .$

\vspace{3mm}

\textbf{Result 3.7. }From (3.1) it follows that 
\begin{equation*}
\left\vert \lambda \right\vert \left\Vert u\right\Vert _{X}\leq C\left\Vert
f\right\Vert _{X}\eqno(3.5)
\end{equation*}%
for all $\lambda \in S_{\varphi _{2}}$ (for some $\varphi _{2}\in \left[
0,\right. \left. \pi \right) $)$.$ Moreover, it is also easy to deduce from
Lemma 3.2 and Lemma 3.3 that%
\begin{equation*}
\left\Vert u\right\Vert _{X}\leq C\left\Vert f\right\Vert _{X}.\eqno(3.6)
\end{equation*}%
Thus, combining (3.5) and (3.6) we obtain the following resolvent estimate 
\begin{equation*}
\left( 1+\left\vert \lambda \right\vert \right) \left\Vert \left( L+\lambda
\right) ^{-1}\right\Vert _{B(X)}\leq C.
\end{equation*}

\vspace{3mm}

\textbf{Remark 3.8. }It is well known that 
\begin{equation*}
W_{p}^{l}\left( R;E\left( A\right) ,E\right) \hookrightarrow L_{p}\left(
R;E\right)
\end{equation*}%
\textbf{\ }whenever 
\begin{equation*}
E\left( A\right) \hookrightarrow E.
\end{equation*}%
Hence, from Result 3.7 it follows that $L$ is a $\varphi $-positive operator
in $L_{p}\left( R;E\right) .$ \ 

As an application of Theorem 2.3 consider the BVP for integro-differential
equations%
\begin{equation*}
\sum\limits_{k=0}^{l}a_{k}\ast \frac{\partial ^{k}u}{\partial x^{k}}%
+\sum\limits_{k=0}^{l}b_{k}\frac{\partial ^{k}u}{\partial x^{k}}%
+\dint\limits_{-\infty }^{\infty }\mu (x-z)\frac{\partial ^{2}u}{\partial
y^{2}}(z,y)dt=f(x,y),\text{ }\eqno(3.7)
\end{equation*}%
\begin{equation*}
u(x,0)=u(x,1)\text{ and }u_{y}(x,0)=u_{y}(x,1),\text{ }x\in (-\infty ,\infty
),\text{ }y\in (0,1)
\end{equation*}%
in a mixed norm space\ $X=L_{p,q}\left( R\times \lbrack 0,1]\right) $ where $%
1<p,q<\infty .$\ It is well-known that the differential expression%
\begin{equation*}
Au=-u^{\prime \prime }+bu
\end{equation*}%
defines a positive operator $A$ acting on $L_{q}\left( [0,1]\right) $ with
domain and satisfying $u(0)=u(1)$ and $u^{\prime }(0)=u^{\prime }(1).$ Thus
by Theorem 3.5 the BVP (3.7) has a unique solution $u\in $ $%
W_{p,q}^{(l,2)}\left( R\times \lbrack 0,1]\right) $ (anisotropic Sobolev
spaces ) and the following coercive inequality is valid:

\begin{equation*}
\sum\limits_{k=0}^{l}\left( \left\Vert a_{k}\ast \frac{\partial ^{k}u}{%
\partial x^{k}}\right\Vert _{X}+\left\Vert \frac{\partial ^{k}u}{\partial
x^{k}}\right\Vert _{X}\right) +\left\Vert \mu \ast \frac{\partial ^{2}u}{%
\partial y^{2}}\right\Vert _{X}\leq C\left\Vert f\right\Vert _{X}.
\end{equation*}

\vspace{3mm}

\section*{4. Application to nonlinear equations}

Here we will apply the main results of previous section to Cauchy problem
for semilinear convolution equation (1.1) and the BVP for nonlinear elliptic
equation (1.2). The main tool we use here is abstract existence and
uniqueness theorem of Amann [15] and Shakhmurov [16]. For the exposition of
abstract quasilinear equations see e.g. [15] and the references therein.

First we show that $L$ is $R$-positive operator.

\vspace{3mm}

\textbf{Theorem 4.0. }Let $E$ be a $UMD$ space and $A$ be a $R$-positive
operator in $E.$ If Condition 3.1 holds and $1<p<\infty $ then $L$ is $R$%
-positive operator in $L_{p}\left( R;E\right) .$

\textbf{Proof. }The Remark 3.8 ensures us $\varphi $-positivity of $L$ in $%
L_{p}\left( R;E\right) $. Therefore it suffices to to show that the
following set 
\begin{equation*}
S=\left\{ (1+\lambda )\left( L+\lambda \right) ^{-1}:\lambda \in S_{\varphi
}\right\}
\end{equation*}%
is $R$-bounded. From the proof of Theorem 3.5 and Result 3.6 we know 
\begin{equation*}
(1+\lambda )\left( L+\lambda \right) ^{-1}f=F^{-1}\left[ \sigma \left( \xi
,\lambda \right) \hat{f}\right] \text{, \ }
\end{equation*}%
for all $f\in L_{p}\left( R;E\right) ,$ where 
\begin{equation*}
\sigma \left( \xi ,\lambda \right) =(1+\lambda )(\hat{\mu}(\xi )+\nu )^{-1}%
\left[ A+\eta \left( \xi \right) +\lambda )\right] ^{-1}.
\end{equation*}%
By virtue of Lemma 3.2 and Lemma 3.3 one can easily prove 
\begin{equation*}
\left\{ \sigma \left( \xi ,\lambda \right) :\xi \in R\backslash \left\{
0\right\} \right\} \text{ and }\left\{ \xi \sigma \left( \xi ,\lambda
\right) :\xi \in \dot{R}\right\}
\end{equation*}%
are $R$-bounded sets and that $\sigma \left( \xi ,\lambda \right) $ is UFM
in $X=L_{p}\left( R;E\right) .$ Thus 
\begin{equation*}
\int\limits_{0}^{1}\left\Vert \sum\limits_{j=1}^{m}r_{j}\left( y\right)
(1+\lambda _{j})\left( L+\lambda _{j}\right) ^{-1}f_{j}\right\Vert
_{X}dy=\int\limits_{0}^{1}\left\Vert \sum\limits_{j=1}^{m}r_{j}\left(
y\right) F^{-1}\left[ \sigma \left( \xi ,\lambda _{j}\right) \hat{f}\right] 
\hat{f}_{j}\right\Vert _{X}dy
\end{equation*}

\begin{equation*}
=\int\limits_{0}^{1}\left\Vert F^{-1}\left[ \sigma \left( \xi ,\lambda
_{j}\right) \sum\limits_{j=1}^{m}r_{j}\left( y\right) \hat{f}\right]
\right\Vert _{X}dy\leq C\int\limits_{0}^{1}\left\Vert
\sum\limits_{j=1}^{m}r_{j}\left( y\right) f_{j}\right\Vert _{X}dy,
\end{equation*}%
for each $m\in \mathbf{N,}$ $(1+\lambda _{j})\left( L+\lambda _{j}\right)
^{-1}\in S,$ $f_{j}\in L_{p}\left( R;E\right) $ and for all independent,
symmetric, $\left\{ -1,1\right\} $-valued random variables $r_{j}$ on $\left[
0,1\right] $. Hence we get the assertion.

Let $E$ be a Banach space, $1<p,q<\infty $ and $\mathbf{p}=(p,q)$. In what
follows, $L_{\mathbf{p}}\left( [0,T]\times R;E\right) $ will denote the
space of all $\mathbf{p}$-summable $E$-valued\ functions with finite mixed
norm 
\begin{equation*}
\left\Vert f\right\Vert _{L_{\mathbf{p}}\left( [0,T]\times R;E\right)
}=\left( \int\limits_{0}^{T}\left( \int\limits_{-\infty }^{\infty
}\left\Vert f\left( t,x\right) \right\Vert _{E}^{q}dx\right) ^{\frac{p}{q}%
}dt\right) ^{\frac{1}{p}}<\infty .
\end{equation*}%
Furthermore, $W_{\mathbf{p}}^{(1,l)}(J\times R;E(A),E)$ denotes anisotropic
Sobolev spaces i.e. the space of all functions $u\in L_{\mathbf{p}}\left(
J\times R;E(A)\right) $ such that $\frac{\partial u}{\partial t}$ and $\frac{%
\partial ^{k}u}{\partial x^{k}}\in L_{\mathbf{p}}\left( J\times R;E\right) ,$
$k=0,\cdot \cdot \cdot ,l$ and 
\begin{equation*}
\left\Vert \frac{\partial u}{\partial t}\right\Vert _{L_{\mathbf{p}}\left(
J\times R;E\right) }+\sum\limits_{k=0}^{l}\left\Vert \frac{\partial ^{k}u}{%
\partial x^{k}}\right\Vert _{L_{\mathbf{p}}\left( J\times R;E\right)
}+\left\Vert Au\right\Vert _{L_{\mathbf{p}}\left( J\times R;E\right)
}<\infty .
\end{equation*}%
It is well known that the Besov spaces has significant embedding properties.
Let us recall some of them: 
\begin{equation*}
W_{q}^{l+1}(X)\hookrightarrow B_{q,r}^{s}(X)\hookrightarrow
W_{q}^{l}(X)\hookrightarrow L_{q}(X)\text{ where }l<s<l+1,
\end{equation*}%
\begin{equation*}
B_{\infty ,1}^{s}(X)\hookrightarrow C^{s}(X)\hookrightarrow B_{\infty
,\infty }^{s}(X)\text{ for }s\in \mathbf{Z},
\end{equation*}%
and%
\begin{equation*}
B_{p,1}^{\frac{N}{p}}(R^{N},X)\hookrightarrow L_{\infty }(R^{N},X)\text{ for 
}s\in \mathbf{Z}.
\end{equation*}%
For the definition of anisotropic Besov spaces see e.g. [23] and [8].

Now let us establish a lemma for representation of iterated spaces, that
will be useful in the proof of our main theorems and their applications. 
\vspace{3mm}

\textbf{Proposition 4.1. }Suppose $\Omega _{1},\Omega _{2}\subseteq R^{n},$ $%
s\in R$ and $1\leq p,q,r\leq \infty .$ Then: 
\begin{equation*}
\text{ }
\end{equation*}%
\textbf{\ }(\textbf{i}) 
\begin{equation*}
B_{p,r}^{s}\left( \Omega _{1},L_{q}(\Omega _{2})\right) =B_{\mathbf{p}%
,r}^{(s,0)}\left( \Omega _{1}\times \Omega _{2}\right) ,
\end{equation*}%
(\textbf{ii})%
\begin{equation*}
L_{p}\left( \Omega _{1},B_{q,r}^{s}(\Omega _{2})\right) =B_{\mathbf{p}%
,r}^{(0,s)}\left( \Omega _{1}\times \Omega _{2}\right)
\end{equation*}%
and%
\begin{equation*}
W_{p}^{l}\left( \Omega _{1},B_{q,r}^{s}(\Omega _{2})\right) =B_{\mathbf{p}%
,r}^{(l,s)}\left( \Omega _{1}\times \Omega _{2}\right)
\end{equation*}%
(\textbf{iii})%
\begin{equation*}
B_{p,r}^{s}\left( \Omega _{1},B_{q,r}^{\sigma }(\Omega _{2})\right) =B_{%
\mathbf{p},r}^{\left( s,\sigma \right) }\left( \Omega _{1}\times \Omega
_{2}\right) .
\end{equation*}

\textbf{Proof. (i) }Taking into consideration [23] and the interpolation
definition of Besov spaces i.e. 
\begin{equation*}
B_{q,r}^{s}(R^{N};X)=\left( L_{q}(R^{N};X),W_{q}^{m}(R^{N};X)\right) _{\frac{%
s}{m},r},
\end{equation*}%
we get the first assertion:%
\begin{equation*}
\begin{array}{lll}
B_{p,r}^{s}\left( \Omega _{1},L_{q}(\Omega _{2})\right) & = & \displaystyle%
\left( L_{p}(\Omega _{1},L_{q}(\Omega _{2})),W_{p}^{m}\left( \Omega
_{1},L_{q}(\Omega _{2})\right) \right) _{\frac{s}{m},r} \\ 
&  &  \\ 
& = & \displaystyle\left( L_{\mathbf{p}}(\Omega _{1}\times \Omega _{2}),W_{%
\mathbf{p}}^{(m,0)}\left( \Omega _{1}\times \Omega _{2}\right) \right) _{%
\frac{s}{m},r} \\ 
&  &  \\ 
& = & \displaystyle B_{\mathbf{p},r}^{(s,0)}\left( \Omega _{1}\times \Omega
_{2}\right) .%
\end{array}%
\end{equation*}%
(\textbf{ii}) Next we apply [22, Theorem 5.1.2] along with definition of
Besov spaces, [23] and (i) we get desired result:%
\begin{equation*}
\begin{array}{lll}
L_{p}\left( \Omega _{1},B_{q,r}^{s}(\Omega _{2})\right) & = & \displaystyle%
\left( L_{p}\left( \Omega _{1},L_{q}(\Omega _{2})\right) ,L_{p}\left( \Omega
_{1},W_{q}^{m}(\Omega _{2})\right) \right) _{\frac{s}{m},r} \\ 
&  &  \\ 
& = & \displaystyle\left( L_{\mathbf{p}}(\Omega _{1}\times \Omega _{2}),W_{%
\mathbf{p}}^{(0,m)}\left( \Omega _{1}\times \Omega _{2}\right) \right) _{%
\frac{s}{m},r} \\ 
&  &  \\ 
& = & \displaystyle B_{\mathbf{p},r}^{(0,s)}\left( \Omega _{1}\times \Omega
_{2}\right) ,%
\end{array}%
\end{equation*}%
and%
\begin{equation*}
\begin{array}{lll}
W_{p}^{l}\left( \Omega _{1},B_{q,r}^{s}(\Omega _{2})\right) & = & %
\displaystyle\left( W_{p}^{l}\left( \Omega _{1},L_{q}(\Omega _{2})\right)
,W_{p}^{l}\left( \Omega _{1},W_{q}^{m}(\Omega _{2})\right) \right) _{\frac{s%
}{m},r} \\ 
&  &  \\ 
& = & \displaystyle\left( W_{\mathbf{p}}^{(l,0)}(\Omega _{1}\times \Omega
_{2}),W_{\mathbf{p}}^{(l,m)}\left( \Omega _{1}\times \Omega _{2}\right)
\right) _{\frac{s}{m},r} \\ 
&  &  \\ 
& = & \displaystyle B_{\mathbf{p},r}^{(l,s)}\left( \Omega _{1}\times \Omega
_{2}\right) .%
\end{array}%
\end{equation*}%
(\textbf{iii}) Finally, with the help of (ii) and [23] we obtain 
\begin{equation*}
\begin{array}{lll}
B_{p,r}^{s}\left( \Omega _{1},B_{q,r}^{\sigma }(\Omega _{2})\right) & = & %
\displaystyle\left( L_{p}\left( \Omega _{1},B_{q,r}^{\sigma }(\Omega
_{2})\right) ,W_{p}^{m}\left( \Omega _{1},B_{q,r}^{\sigma }(\Omega
_{2})\right) \right) _{\frac{s}{m},r} \\ 
&  &  \\ 
& = & \displaystyle\left( B_{\mathbf{p},r}^{(0,\sigma )}(\Omega _{1}\times
\Omega _{2}),B_{\mathbf{p,r}}^{(m,\sigma )}\left( \Omega _{1}\times \Omega
_{2}\right) \right) _{\frac{s}{m},r} \\ 
&  &  \\ 
& = & \displaystyle B_{\mathbf{p},r}^{(s,\sigma )}\left( \Omega _{1}\times
\Omega _{2}\right) .%
\end{array}%
\end{equation*}

Let $X$ and $Y$ be nonempty sets and $J_{T}=[0,T)$. A nonlocal map $%
F:X^{J}\rightarrow Y^{J}$ is said to have Voltera property (or Voltera map)
if for every $S\in \mathring{J}_{T}$ and $u\in $ $X$%
\begin{equation*}
F(u)|_{J_{S}}=F(u|_{J_{S}}).
\end{equation*}%
It is clear that local maps satisfy Voltera property automatically.

By a solution on $J_{S}$ of (1.1) we mean $u\in W_{\mathbf{p}%
,loc}^{(1,l)}(J\times R;E(A),E)$ such that $u|_{J_{S}}$ belongs to $W_{%
\mathbf{p}}^{(1,l)}(J_{S}\times R;E(A),E)$ for each $S\in \mathring{J}_{T}$
and satisfies a.e. the semilinear convolution operator equation (1.1).

It is known that if $u$ $\in W_{\mathbf{p}}^{(1,l)}(J_{S}\times R;E(A),E)$
then 
\begin{equation*}
\frac{\partial ^{j}u}{\partial x^{j}}\in W_{\mathbf{p}}^{(1,l-j)}\left(
J\times R;\left( E(A),E\right) _{\frac{j}{l}}\right) ,\text{ }1\leq j\leq
l-1.
\end{equation*}%
Let $B_{0}$ be a cartesian product of above spaces, namely%
\begin{equation*}
B_{0}=\prod\limits_{j=0}^{l-1}W_{\mathbf{p}}^{(1,l-j)}\left( J\times
R;\left( E(A),E\right) _{\frac{j}{l}}\right) .
\end{equation*}%
\bigskip

\textbf{Theorem 4.2. }Let $E$ be a $UMD$ space, $A$ be a $R$-positive
operator in $E$ with $\varphi \in \left( \frac{\pi }{2},\pi \right) $ and $%
J=[0,T).$ Assume Condition 3.1 holds, $1<p,q<\infty $ and $u^{0}\in
B_{q,p}^{l/p^{\prime }}\left( R;E\right) \cap L_{q}\left( R;\left(
E,E(A)\right) _{1/p^{\prime },p}\right) $. If

(i) $F$ is a Voltera map from $B_{0}$ to $L_{\mathbf{p}}(J\times R,E),$ and

(ii) there exists $r\in (p,\infty ]$ such that $F-F(0)$ is uniformly
Lipschitz continuous on bounded subsets of $B_{0}$ with values in 
\begin{equation*}
L_{\mathbf{r}}(J\times R,E),
\end{equation*}%
then $\left( 1.1\right) $ has a unique maximal solution. The maximal
interval of existence, $J_{\max }$ is open in $J.$ Moreover, if for the
unique maximal solution $u$ of (1.1) $F(u)\in L_{\mathbf{p}}(J_{\max }\times
R,E),$ then $J_{\max }=J.$

\textbf{Proof.} Let $E_{1}=W_{q}^{l}\left( R;E\left( A\right) ,E\right) $
and $E_{0}=L_{q}\left( R;E\right) $ where 
\begin{equation*}
E\left( A\right) \hookrightarrow E.
\end{equation*}%
\ By the intersection property of interpolation, 
\begin{equation*}
(E_{0},E_{1})_{1/p^{\prime },p}=B_{q,p}^{l/p^{\prime }}\left( R;E\right)
\cap L_{q}\left( R;\left( E,E(A)\right) _{1/p^{\prime },p}\right) .
\end{equation*}%
It is clear that 
\begin{equation*}
W_{q}^{l}\left( R;E\left( A\right) ,E\right) \hookrightarrow
(E_{0},E_{1})_{1/p^{\prime },p}\hookrightarrow L_{q}\left( R;E\right)
\end{equation*}%
and 
\begin{equation*}
W_{p}^{1}(J;E_{1},E_{0})=W_{\mathbf{p}}^{(1,0)}(J\times R,E)\cap W_{\mathbf{p%
}}^{(0,l)}(J\times R,E(A),E)=W_{\mathbf{p}}^{(1,l)}(J\times R;E(A),E).
\end{equation*}%
Now consider the linear form of (1.1), 
\begin{equation*}
\frac{\partial u}{\partial t}+\sum\limits_{k=0}^{l}a_{k}\ast \frac{\partial
^{k}u}{\partial x^{k}}+\sum\limits_{k=0}^{l}b_{k}\frac{\partial ^{k}u}{%
\partial x^{k}}+\mu \ast Au+\nu Au=f(t,x),\eqno(4.1)
\end{equation*}%
\begin{equation*}
u(0,x)=u_{0}\text{, }t\in \left( 0,T\right) ,\text{ }x\in \left( -\infty
,\infty \right)
\end{equation*}%
in $L_{\mathbf{p}}\left( J\times R;E\right) =L_{p}(J;E_{0}).$ By the Theorem
4.1 $L$ is $R$-positive operator in $L_{q}\left( R;E\right) $ with $\varphi
\in \left( \frac{\pi }{2},\pi \right) $ and thus a generator of analytic
semigroup. Since $\left( 4.1\right) $\ can be written in an abstract form 
\begin{equation*}
\frac{du}{dt}+Lu=f(t,x),\text{ }u\left( 0\right) =u^{0},\text{ }t\in J,\eqno%
(4.2)
\end{equation*}%
$\left[ \text{20, Theorem 4.2}\right] $ and Theorem 3.5 ensures that, for
each $f\in Y=L_{\mathbf{p}}\left( J\times R;E\right) ,$ $\left( 4.1\right) $
has a unique solution $u\in W_{\mathbf{p}}^{(1,l)}(J\times R;E(A),E)$ and
the following coercive estimate%
\begin{equation*}
\left\Vert \frac{\partial u}{\partial t}\right\Vert
_{Y}+\sum\limits_{k=0}^{l}\left( \left\Vert a_{k}\ast \frac{\partial ^{k}u}{%
\partial x^{k}}\right\Vert _{Y}+\left\Vert \frac{\partial ^{k}u}{\partial
x^{k}}\right\Vert _{Y}\right) +\left\Vert \mu \ast Au\right\Vert
_{Y}+\left\Vert Au\right\Vert _{Y}\leq C\left( \left\Vert f\right\Vert
_{Y}+\left\Vert Au_{0}\right\Vert _{E}\right)
\end{equation*}%
holds. It is known that if $A$ is a constant map, the maximal regularity of
problem (4.1) is independent of bounded intervals $J$ and of $p$ (see $\left[
\text{20, Remarks 6.1(d) and (e)}\right] $). In a similar manner nonlinear
equation (1.1) reduces to the following abstract form%
\begin{equation*}
\frac{du}{dt}+Lu=\bar{F}(u),\text{ }u\left( 0\right) =u^{0},\text{ }t\in J.
\end{equation*}%
From assumptions (i) and (ii) it follows that $\bar{F}$ is a Voltera map
from $W_{p}^{1}(J;E_{1},E_{0})$ to $L_{p}(J;E_{0}),$ and $F-F(0)$ is
uniformly Lipschitz continuous on bounded subsets of $%
W_{p}^{1}(J;E_{1},E_{0})$ with values in $L_{r}(J;E_{0}).$ Taking into
account the fact that $L$ has maximal regularity property independent of
intervals $J,$ by [15, Theorem 2.1], $\left( 1.1\right) $ has a unique
maximal solution $u$ and this solution is globally defined whenever $u\in
W_{p}^{1}(J_{\max }).$ \hbox{\vrule height7pt width5pt}

Now let us study (4.1) and (1.1) in concrete settings.\vspace{3mm}

\textbf{Example 4.3. }Suppose $A$ is an $n\times n$ matrix whose eigenvalues
have positive real parts. Let $a=e^{-k\left\vert x\right\vert },$ $k>0,$ $%
b\in R$ and $1<p,q<\infty .$\textbf{\ }It is clear that 
\begin{equation*}
F(e^{-k\left\vert x\right\vert })=\hat{a}(\xi )=\frac{2i\xi }{k^{2}+\xi ^{2}}
\end{equation*}%
and 
\begin{equation*}
\eta \left( \xi \right) =-\xi ^{2}(\hat{a}(\xi )+b)\in S_{\varphi },\text{ }%
\varphi \in \left[ 0,\right. \left. \pi \right) .
\end{equation*}%
Moreover, $\hat{a}(\xi )$ satisfy all assumptions of Condition 3.1. Now
consider the system of integro-differential equations 
\begin{equation*}
u_{t}+a\ast u_{xx}+bu_{xx}+Au=f(t,x),\text{ \ }t\in (0,T),\text{ }-\infty
<x<\infty \eqno(4.3)
\end{equation*}%
\begin{equation*}
u(0,x)=u_{0}\in B_{q,p}^{l/p^{\prime }}\left( R,R^{n}\right) \text{, }x\in
\left( -\infty ,\infty \right) .
\end{equation*}%
Since all assumptions of the Theorem 4.2 are satisfied, for each $f\in Y=L_{%
\mathbf{p}}\left( J\times R,R^{n}\right) ,$ (4.3) has a unique solution $%
u\in W_{\mathbf{p}}^{(1,2)}\left( J\times R,R^{n}\right) $ and the coercive
inequality holds%
\begin{equation*}
\left\Vert u_{t}\right\Vert _{Y}+\left\Vert a\ast u_{xx}\right\Vert
_{Y}+\left\Vert u_{xx}\right\Vert _{Y}+\left\Vert u\right\Vert _{Y}\leq
C\left( \left\Vert f\right\Vert _{Y}+\left\Vert u_{0}\right\Vert
_{B_{q,p}^{l/p^{\prime }}}\right) .
\end{equation*}

\vspace{3mm}

\textbf{Example 4.4. }Assume $a_{k},b_{k},\mu $ and $\nu $ satisfy
assumptions of Condition 3.1, $l>2$ and $f$ is a nonlocal map as in the
Theorem 4.2. Let $A=-\Delta +c$ and $E=L_{q}(\Omega )$\textbf{\ }where $%
\Omega \subset R^{2}$ has sufficiently smooth boundary.\textbf{\ }From
Proposition 4.1 we have%
\begin{eqnarray*}
B_{q,p}^{l/p^{\prime }}\left( R,L_{q}(\Omega )\right) \cap L_{q}\left( R, 
\left[ L_{q}\left( \Omega \right) ,W_{q}^{2}\left( \Omega \right) \right]
_{1/p^{\prime },p}\right) &=&B_{q,p}^{l/p^{\prime }}\left( R\times \Omega
\right) \cap B_{q,p}^{2/p^{\prime }}\left( R\times \Omega \right) \\
&=&B_{q,p}^{l/p^{\prime }}\left( R\times \Omega \right) .
\end{eqnarray*}%
Now consider an initial-boundary value problem for 3D semilinear
integro-differential equation 
\begin{equation*}
\frac{\partial u}{\partial t}+\sum\limits_{k=0}^{l}\left( a_{k}\ast \frac{%
\partial ^{k}u}{\partial x^{k}}+b_{k}\frac{\partial ^{k}u}{\partial x^{k}}%
\right) -\mu \ast \Delta u-\nu \Delta u=f\left( u,u_{x},\cdot \cdot \cdot ,%
\frac{\partial ^{l-1}u}{\partial x^{l-1}}\right) ,\text{ \ }t\in (0,T)\eqno%
(4.4)
\end{equation*}%
\begin{equation*}
u(0,x,y,z)=u_{0}\text{, }x\in \left( -\infty ,\infty \right) \text{, }%
(y,z)\in \Omega
\end{equation*}%
\begin{equation*}
u(t,x,y,z)|_{\partial \Omega }=0.
\end{equation*}%
Therefore by Theorem 4.2, if $u_{0}\in B_{q,p}^{l/p^{\prime }}\left( R\times
\Omega \right) $ then (4.4) has a unique maximal solution $u\in W_{\mathbf{p}%
}^{\mathbf{s}}\left( J\times R\times \Omega \right) $ where $\mathbf{s}%
=(1,l,2)$ and $\mathbf{p}=(p,q).$ \hbox{\vrule
height7pt width5pt}

Let $1<p,q<\infty ,$ $J_{0}=[0,T_{0}]$ and $J=[0,T].$ Suppose boundary
condition of (1.2) is nondegenerate i.e.%
\begin{equation*}
\left\vert 
\begin{array}{cc}
\alpha _{1} & \alpha _{2} \\ 
\beta _{1} & \beta _{2}%
\end{array}%
\right\vert \neq 0,
\end{equation*}%
$E_{1}=W_{q}^{l}\left( R;E\left( A\right) ,E\right) ,$ $E_{0}=L_{q}\left(
R;E\right) $ and $Y=W_{p}^{2}\left( \left( 0,T\right) ;E_{1},E_{0}\right) .$
From the well known trace theorem [24] we know that if 
\begin{equation*}
u\in W_{p}^{2}\left( \left( 0,T\right) ;E_{1},E_{0}\right)
\end{equation*}%
then%
\begin{equation*}
u(0)\in X_{0}=\left( E_{1},E_{0}\right) _{\frac{1}{2p},p}=B_{q,p}^{l\frac{%
\left( 2p-1\right) }{2p}}\left( R;E\right) \cap L_{q}\left( R;\left(
E(A),E\right) _{\frac{1}{2p},p}\right)
\end{equation*}%
and%
\begin{equation*}
u^{\prime }(0)\in X_{1}=\left( E_{1},E_{0}\right) _{\frac{p+1}{2p}%
,p}=B_{q,p}^{l\frac{\left( p-1\right) }{2p}}\left( R;E\right) \cap
L_{q}\left( R;\left( E(A),E\right) _{\frac{p+1}{2p},p}\right) .
\end{equation*}%
Let $\bar{u}$ $=\left( u,\frac{\partial u}{\partial t}\right) $ and $\Omega
_{0}$, $\Omega $ denote infinite strips $[0,T)\times R$ and $[0,T_{0})\times
R$ respectively and $B_{0}$ be a cartesian product of above spaces i.e.%
\begin{equation*}
B_{0}=X_{0}\times X_{1}.
\end{equation*}

\vspace{3mm}

\textbf{Theorem 4.5. }Let $E$ be a $UMD$ space and $A$ be a $R$-positive
operator in $E.$ Assume Condition 3.1 holds and $F:\Omega _{0}\times
B_{0}\rightarrow L_{q}\left( R;E\right) $ satisfies Caratheodory and
Lipschitz continuity conditions i.e.

(i) $F(\cdot ,\bar{u})$ is measurable for each $\bar{u}\in B_{0}$, $%
F(t,\cdot )$ is continuous almost all $t\in J_{0}=[0,T_{0}),$

(ii) $f(t)=F(t,0)\in L_{\mathbf{p}}(\Omega _{0},E),$ for each $R>0$ there is
a function $\phi _{R}\in L_{p}(J_{0})$ such that 
\begin{equation*}
\left\vert F(t,\bar{u})-F(t,\bar{v})\right\vert _{L_{q}\left( R;E\right)
}\leq \phi _{R}(t)\left( \left\vert \bar{u}-\bar{v}\right\vert
_{B_{0}}\right) ,\text{ a.a }t\in J,\text{ }\bar{u},\bar{v}\in
B_{0},\left\vert \bar{u}\right\vert _{B_{0}},\left\vert \bar{v}\right\vert
_{B_{0}}\leq R.
\end{equation*}%
If $f_{1},f_{2}\in X_{1},$ then there exist $T\in J_{0}$ such that $\left(
1.2\right) $ admits a unique solution $u\in W_{\mathbf{p}}^{(2,l)}(\Omega
;E(A),E).$

\textbf{Proof.} First we study linear form of (1.2) in the moving boundary $%
(0,b(\xi ))$ i.e. 
\begin{equation*}
-\frac{\partial ^{2}u}{\partial t^{2}}+\sum\limits_{k=0}^{l}a_{k}\ast \frac{%
\partial ^{k}u}{\partial x^{k}}+\sum\limits_{k=0}^{l}b_{k}\frac{\partial
^{k}u}{\partial x^{k}}+\mu \ast Au+\nu Au=f(t,x),\eqno(4.5)
\end{equation*}%
\begin{equation*}
\alpha _{1}u(0,x)+\beta _{1}\frac{\partial u(0,x)}{\partial t}=f_{1}(x)\text{%
, }\alpha _{2}u(b(\xi ),x)+\beta _{2}\frac{\partial u(b(\xi ),x)}{\partial t}%
=f_{2}(x),\text{ }x\in \left( -\infty ,\infty \right) .
\end{equation*}%
By the Theorem 4.1, $L$ is $R$-positive operator in $L_{q}\left( R;E\right)
. $ Since $\left( 4.5\right) $\ can be written in abstract form 
\begin{equation*}
-u^{\prime \prime }(t)+Lu\left( t\right) =f\left( t\right)
\end{equation*}%
\begin{equation*}
\alpha _{1}u(0)+\beta _{1}u^{\prime }(0)=f_{1}\text{, }\alpha _{1}u(b(\xi
))+\beta _{2}u^{\prime }(b(\xi ))=f_{2},
\end{equation*}%
$\left[ \text{24, Theorem 4.1}\right] $ and Theorem 3.5 ensures us that for
each $f\in Y_{b}=L_{\mathbf{p}}\left( (0,b(\xi ))\times R;E\right) ,$ $%
\left( 4.1\right) $ has a unique solution $u\in W_{\mathbf{p}%
}^{(2,l)}((0,b(\xi ))\times R;E(A),E)$. Similarly, nonlinear BVP (1.2)\
reduces to abstract form 
\begin{equation*}
-u^{\prime \prime }(t)+Lu\left( t\right) =F\left( t,u,u^{\prime }\right)
\end{equation*}%
\begin{equation*}
\alpha _{1}u(0)+\beta _{1}u^{\prime }(0)=f_{1}\text{, }\alpha _{1}u(T)+\beta
_{2}u^{\prime }(T)=f_{2}.
\end{equation*}%
Since $L$ is $R$-positive operator in $E_{0}=L_{q}\left( R;E\right) ,$ by $%
\left[ \text{16, Theorem 4}\right] $ there exist $T\in \lbrack 0,T_{0})$
such that $\left( 1.2\right) $ admits a unique solution $u\in
W_{p}^{2}(J;E_{1},E_{0})=W_{\mathbf{p}}^{(2,l)}(\Omega ;E(A),E).$%
\hbox{\vrule
height7pt width5pt}

Finally, choosing $p=q=2$ and $l=4$ in (1.2) we consider BVP for the system
of nonlinear equations:%
\begin{equation*}
-\frac{\partial ^{2}u}{\partial t^{2}}+a\ast \frac{\partial ^{4}u}{\partial
x^{4}}+b\frac{\partial ^{4}u}{\partial x^{4}}+Au=f\left( u,\nabla u,\frac{%
\partial ^{2}u}{\partial x^{2}},\frac{\partial ^{3}u}{\partial x^{3}}\right)
,\text{ }t\in (0,T_{0}),\text{ }-\infty <x<\infty \eqno(4.6)
\end{equation*}%
\begin{equation*}
\alpha _{1}u(0,x)+\beta _{1}\frac{\partial u(0,x)}{\partial t}=f_{1}(x)\text{%
, }\alpha _{2}u(T,x)+\beta _{2}\frac{\partial u(T,x)}{\partial t}=f_{2}(x),%
\text{ }x\in \left( -\infty ,\infty \right) ,
\end{equation*}%
where $A,$ $a$ and $b$ are as in Example 4.3.\textbf{\ }If $f_{1}$\textbf{, }%
$f_{2}\in B_{2,2}^{3}(R,R^{n})$\textbf{\ }then\textbf{\ }by Theorem 4.5,
there exist $T\in \left( 0,T_{0}\right) $ such that (4.6) admits a unique
solution $u\in W_{2}^{(2,4)}\left( [0,T]\times R,R^{n}\right) .$

There are a lot of positive operators in concrete Banach spaces. Therefore,
choosing concrete positive differential, pseudo differential operators, or
finite, infinite matrices, etc.\ instead of $A,$\ we can obtain global
existence results for various nonlinear convolution equations.

\end{document}